\documentclass[12pt]{article}
\usepackage{epic,eepic,epsf,epsfig}
\usepackage{amsmath}
\usepackage{amssymb}
\usepackage{amsbsy}

\usepackage{amsfonts}%,srcltx}
\newtheorem{theorem}{Theorem}[section]%

\def\f{\noindent}

\newcommand{\qed}{\mbox{\raisebox{0.7ex}{\fbox{}}} \vspace{4truemm}}
\def\demo{\f {\bf Proof.}\hskip10pt}

\begin{document}

\baselineskip 16pt

\title{ \vspace{-1.2cm}
On an extension of Shlyk's theorem
\thanks{\scriptsize This work was supported by Shandong Provincial Natural Science Foundation, China (ZR2017MA022) and NSFC (11761079).
\newline
 \hspace*{0.5cm} \scriptsize{E-mail addresses:} jtshi2005@163.com\,(J. Shi),\,\,xufj2023@s.ytu.edu.cn\,(F. Xu),\,\,ln18865550588@163.com\,(N. Li).}}

\author{Jiangtao Shi$\,^a$,\,\,Fanjie Xu$\,^a$,\,\,Na Li$\,^b$\\
\\
{\small{\em $^a$ School of Mathematics and Information Sciences,
Yantai University, Yantai 264005, P.R. China}}\\
{\small{\em $^b$ College of Mathematics and Statistics,  Zaozhuang University, Zaozhuang 277160, P.R. China}}}

\date{ }

\maketitle \vspace{-.8cm}

\begin{abstract}
In this paper, we prove that the intersection of all non-nilpotent maximal subgroups of a non-solvable group containing the normalizer of some Sylow subgroup is nilpotent,
which provides an extension of Shlyk's theorem.

\medskip
\f {\bf Keywords:} non-solvable group; non-nilpotent maximal subgroup; normalizer; Sylow subgroup; nilpotent
\end{abstract}

\section{Introduction}
All groups in this paper are considered to be finite. It is known that the intersection of all maximal subgroups of a group equals its Frattini subgroup. As a generalization, Gasch$\rm\ddot{u}$tz
{\rm\cite{gasch}} showed that the intersection of all non-normal maximal subgroups of a non-nilpotent group is nilpotent. Furthermore, Shidov {\rm\cite{shidov}} proved that the intersection of all
non-nilpotent maximal subgroups of a non-solvable group is nilpotent.

It is easy to see that the intersection of all non-nilpotent maximal subgroups of a non-solvable group also equals its Frattini subgroup by Shidov's theorem. Therefore, the intersection of 
all non-nilpotent maximal subgroups of a non-solvable group $G$ equals the intersection of all maximal subgroups of $G$. However, if a group $G$ is solvable,
we cannot ensure that the intersection of all non-nilpotent maximal subgroups of $G$ equals its Frattini subgroup. For example, let $G=\langle a,\,b\mid\,a^6=b^2=1,\,b^{-1}ab=a^{-1}\rangle$ be a
dihedral group of order 12. It is obvious that the intersection of all non-nilpotent maximal subgroups of $G$ equals $\langle a^2\rangle$ but $\it\Phi$$(G)$=1.

Note that a group in which every maximal subgroup is nilpotent is solvable by {\rm\cite{redei1956}} and a group in which all non-nilpotent maximal
subgroups are normal is also solvable by {\rm\cite{lishi}}, then any non-solvable group must have non-normal non-nilpotent maximal subgroups.
As an extension of Gasch$\rm\ddot{u}$tz's theorem and Shidov's theorem, Shlyk {\rm\cite{shlyk}} had the following result.

\begin{theorem} {\rm\cite{shlyk}}\ \  The intersection of all non-normal non-nilpotent maximal subgroups of a non-solvable group $G$ is nilpotent.
\end{theorem}

Note that {\rm\cite[Lemma 6]{shizhangguo}} showed that the intersection of all non-normal non-nilpotent maximal subgroups of a non-solvable group $G$ equals the intersection of all non-normal
maximal subgroups of $G$.

Although the intersection of all non-normal non-nilpotent maximal subgroups of a non-solvable group is nilpotent by Shlyk's theorem, we cannot
ensure that the intersection of all non-normal non-nilpotent maximal subgroups of a non-solvable group equals its Frattini subgroup. For example, let $G=A_5\times Z_7$. It is obvious that
the intersection of all non-normal non-nilpotent maximal subgroups of $G$ equals $Z_7$ but $\it\Phi$$(G)$=1.

It is easy to see that any maximal subgroup of a group $G$ containing the normalizer of some Sylow subgroup must be a non-normal maximal subgroup of $G$ but a non-normal maximal subgroup of $G$
might not be a maximal subgroup of $G$ containing the normalizer of some Sylow subgroup. In this paper, as a further extension of Shlyk's theorem, we obtain the following result,
the proof of which is given in Section~\ref{s2}.

\begin{theorem}\ \ \label{th2} The intersection of all non-nilpotent maximal subgroups of a non-solvable $G$ containing the normalizer of some Sylow subgroup is nilpotent.
\end{theorem}

In Section~\ref{s3}, we give a simple proof to show that the following result also holds.

\begin{theorem}\ \ \label{th3} The intersection of normalizers of all non-normal Sylow subgroups of a non-nilpotent group $G$ is nilpotent.
\end{theorem}

\section{Proof of Theorem~\ref{th2}}\label{s2}

\demo Since $G$ is non-solvable, $G$ must have non-nilpotent maximal subgroups containing the normalizer of some Sylow subgroup by {\rm\cite[Theorem 1.3]{shixu}}.
Let $N$ be the intersection of all non-nilpotent maximal subgroups of $G$ containing the normalizer of some Sylow subgroup. It is obvious that $N\unlhd G$.

Let $P$ be any Sylow subgroup of $G$. Then $P\cap N$ is a Sylow subgroup of $N$.

Claim that either $P\cap N\unlhd G$ or there exists a nilpotent maximal subgroup $M$ of $G$ such that $N_G(P\cap N)\leq M$.

Otherwise, assume that $P\cap N\ntrianglelefteq G$ and $N_G(P\cap N)$ aren't contained in any nilpotent maximal subgroup of $G$. One has $N_G(P)\leq N_G(P\cap N)<G$. By Frattini argument,
$G=NN_G(P\cap N)$. Let $K$ be a maximal subgroup of $G$ such that $N_G(P\cap N)\leq K$. By the assumption, $K$ is non-nilpotent. Moreover, since $N_G(P)\leq N_G(P\cap N)\leq K$, one has $N\leq K$.
It follows that $G=NN_G(P\cap N)=K$, a contradiction. Therefore, we have
that either $P\cap N\unlhd G$ or there exists a nilpotent maximal subgroup $M$ of $G$ such that $N_G(P\cap N)\leq M$.

Assume $|N|={p_1}^{\alpha_1}{p_2}^{\alpha_2}\cdots{p_s}^{\alpha_s}$, where $p_1,\,p_2,\,\cdots,\,p_s$ are distinct prime divisors of $|N|$, $\alpha_i\geq 1$ for every $1\leq i\leq s$.
Let $P_i\in{\rm Syl}_{p_i}(G)$, $1\leq i\leq s$.

Case (1). If there exists a $P_i$ for some $1\leq i\leq s$ such that $P_i\unlhd G$. Let $T$ be any nilpotent maximal subgroup of $G$. Since $G$ is non-solvable,
one has $P_i\leq T$. It follows that every Sylow subgroup of $N$ being normal in $G$ is contained in any nilpotent maximal subgroup of $G$.

Case (2). If there exists a $P_j$ for some $1\leq j\leq s$ such that $P_j\ntrianglelefteq G$. Then there exists a nilpotent maximal subgroup $M$ of $G$ such that $N_G(P_j)\leq M$, which
implies that $P_j\leq M$.

By cases (1) and (2), one has that every Sylow subgroup of $N$ is contained in some nilpotent maximal subgroup of $G$.

Since $G$ is non-solvable, every nilpotent maximal subgroup of $G$ has even order by {\rm\cite[Theorem 10.4.2]{rob}}. Let $L$ be any nilpotent maximal subgroup of $G$ and $Q$ be Sylow 2-subgroup of $L$.
If $Q$ is not a Sylow 2-subgroup of $G$. Assume that $Q_0$ is a Sylow 2-subgroup of $G$ such that $Q<Q_0$. It follows that $L<N_G(Q)$, which implies that $N_G(Q)=G$. Then $Q\unlhd G$.
One has that $L/Q$ is a nilpotent maximal subgroup of $G/Q$ of odd order. By {\rm\cite[Theorem 10.4.2]{rob}}, $G/Q$ is solvable. Then $G$ is solvable, a contradiction.
Therefore, the Sylow 2-subgroup $Q$ of $L$ is also a Sylow 2-subgroup of $G$. Moreover, $Q\ntrianglelefteq G$. One has $N_G(Q)=L$.

Since all Sylow 2-subgroups of $G$ are conjugate in $G$ by Sylow's theorem, normalizers of all Sylow 2-subgroups of $G$ are also conjugate in $G$.
It follows that all nilpotent maximal subgroups of $G$ are conjugate in $G$.

For every Sylow subgroup $P_i$ of $N$. Note that $N\unlhd G$, one has ${P_i}^g\leq N^g=N$ for any $g\in G$. That is, ${P_i}^g$ is also a Sylow $p_i$-subgroup of $N$ for any $g\in G$.
Since every Sylow subgroup of $N$ is contained in some nilpotent maximal subgroup of $G$ and all nilpotent maximal subgroups of $G$ are conjugate in $G$, there exist
some $P_{10},\,P_{20},\,\cdots,\,P_{s0}$ which are Sylow $p_1$-subgroup, Sylow $p_2$-subgroup, $\cdots$, Sylow $p_s$-subgroup of $N$ respectively and some nilpotent maixmal subgroup $R$ of $G$
such that $P_{i0}\leq R$ for every $1\leq i\leq s$. Since $R$ is nilpotent, one has $[P_{i0},P_{j0}]=1$ for every $1\leq i\neq j\leq s$. It follows that $N$ is nilpotent.\hfill\qed

\section{Proof of Theorem~\ref{th3}}\label{s3}

\demo Since $G$ is non-nilpotent, $G$ must have non-normal Sylow subgroups. Let $N$ be the intersection of normalizers of all non-normal Sylow subgroups of $G$. It is easy to see that $N\unlhd G$.

Let $P$ be any Sylow subgroup of $G$. Then $P\cap N$ is a Sylow subgroup of $N$. Claim that $P\cap N\unlhd G$.

Otherwise, assume $P\cap N\ntrianglelefteq G$. It follows that $P\ntrianglelefteq G$. Then $N_G(P)\leq N_G(P\cap N)<G$. By Frattini argument, one has $G=NN_G(P\cap N)$. Note that $N\leq N_G(P)$ by the assumption.
Then $G=N_G(P)N_G(P\cap N)=N_G(P\cap N)$, this contradicts $P\cap N\ntrianglelefteq G$. Therefore, one has $P\cap N\unlhd G$.

Hence every Sylow subgroup of $N$ is normal in $N$ and then $N$ is nilpotent.\hfill\qed

\end{document}